\documentclass[times,sort&compress,3p,12pt]{elsarticle}
\journal{}

\makeatletter
\def\ps@pprintTitle{%
 \let\@oddhead\@empty
 \let\@evenhead\@empty
 \def\@oddfoot{\centerline{\thepage}}%
 \let\@evenfoot\@oddfoot}
\makeatother

\usepackage{amsmath,amsfonts,amssymb,amsthm,mathtools,hyperref}

\usepackage{kpfonts} 
\usepackage{dsfont} 
\usepackage{xcolor}
\usepackage{geometry}
\newtheorem{proposition}{Proposition}
\newtheorem{corollary}{Corollary}

\newtheorem{remark}{Remark}

\newcommand{\N}{\mathbb{N}}
\newcommand{\R}{\mathbb{R}}
\newcommand{\PP}{\mathsf{P}} 
\newcommand{\QQ}{\mathsf{Q}} 
\newcommand{\EE}{\mathsf{E}} 
\newcommand{\bb}[1]{\boldsymbol{#1}}
\newcommand{\OO}{\mathcal{O}}
\newcommand{\rd}{{\rm d}}
\newcommand{\ind}{\mathds{1}}

\allowdisplaybreaks

\begin{document}

\pagestyle{empty}

\begin{frontmatter}

    \title{Deficiency bounds for the multivariate inverse hypergeometric distribution}%

    \author[a1,a2]{Fr\'ed\'eric Ouimet}%

    \address[a1]{Centre de recherches math\'ematiques, Universit\'e de Montr\'eal, Montr\'eal (Qu\'ebec) Canada H3T 1J4}%
    \address[a2]{Department of Mathematics and Statistics, McGill University, Montr\'eal (Qu\'ebec) Canada H3A 0B9\vspace{-6mm}}%

    \ead{frederic.ouimet2@mcgill.ca}%

    \begin{abstract}
        The multivariate inverse hypergeometric (MIH) distribution is an extension of the negative multinomial (NM) model that accounts for sampling without replacement in a finite population. Even though most studies on longitudinal count data with a specific number of `failures' occur in a finite setting, the NM model is typically chosen over the more accurate MIH model. This raises the question: How much information is lost when inferring with the approximate NM model instead of the true MIH model? The loss is quantified by a measure called deficiency in statistics. In this paper, asymptotic bounds for the deficiencies between MIH and NM experiments are derived, as well as between MIH and the corresponding multivariate normal experiments with the same mean-covariance structure. The findings are supported by a local approximation for the log-ratio of the MIH and NM probability mass functions, and by Hellinger distance bounds.
    \end{abstract}

    \begin{keyword}
        asymptotic statistics \sep comparison of experiments \sep deficiency \sep Hellinger distance \sep local approximation \sep multivariate inverse hypergeometric distribution \sep negative multinomial distribution
        \MSC[2020]{Primary: 62B15, 62E20; Secondary: 60E05, 60F99, 62H10, 62H12}
    \end{keyword}

\end{frontmatter}

\vspace{-4mm}
\section{Introduction}\label{sec:intro}

    The multivariate inverse hypergeometric (MIH) distribution is a flexible distribution introduced by \citet{MR0040607} which has its roots in a variety of practical real-world scenarios involving categorical data. For a finite population size of $N\in \N$ and a given number of $d+1$ disjoint categories with $d\in \N = \{1,2,\ldots\}$, let $p_1,\dots,p_{d+1}\in N^{-1} \, \N$ be a set of probability weights such that $p_{d+1} \equiv q = 1 - (p_1 + \dots + p_d) > 0$. If one assumes that all objects in a given population can be categorized into $d+1$ categories, the $i$-th category having a total of $N p_i$ objects, and one selects a random sample without replacement from the population until $k_{d+1} \equiv n$ objects from the $(d+1)$-th category appear in the sample, then the random vector $\bb{K} = (K_1,\dots,K_d)$ that counts the number of sampled objects that belong to the first $d$ categories follows a MIH distribution. The $(d+1)$-th category is sometimes referred to as the `failure' category and the object placed in it as `failures' in the context where an object fails to be placed in one of the first $d$ categories under study.

    \vspace{2mm}
    When the population under study has infinite size ($N = \infty$), the above count vector $\bb{K} = (K_1,\dots,K_d)$ follows a negative multinomial (NM) distribution. The negative multinomial distribution naturally lends itself to modeling count data in a plethora of domains such as ecology, medicine and manufacturing. For instance, in ecological and biodiversity studies, \citet{Power2016} explains how it could be utilized in the context of linear regression to model the counts of individual encounters of each of $d$ species before a predetermined number of encounters has occurred in another category of species (the `failure' category), effectively aiding in the quantification of biodiversity or species richness. Similarly, in the realm of clinical trials and medical research, \citet{Chiarappa2019phd,MR4242807} have shown that it can serve to model the numbers of adverse events in multiple experimental groups until a set number of adverse events are observed in the control group. In manufacturing, the negative multinomial finds relevance in quality control processes, being used to estimate the number of inspections needed to identify a specific number of defects classified based on varying criteria. Plenty of other examples are investigated in the literature, such as its application to polarimetric image processing \citep{MR3037959}, autoregressive regression models in psychology \citep{doi:10.1080/01621459.1999.10474178,MR1681133}, models of cancer incidence by region \citep{MR1475055}, the analysis of pollen counts \citep{doi:10.2307/2333745}, models of median crossover accident frequencies in transportation research \citep{doi:10.3141/1840-22}, regression models for clustered counts of surgical procedures \citep{doi:10.2307/271020}, the analysis of RNA sequencing \citep{Kusi-Appiah2016phd}, etc. Overall, the negative multinomial distribution provides a flexible and powerful tool for modeling count data in a variety of applications in part because it can also handle overdispersion, zero-inflation and multiple event types.

    \vspace{2mm}
    In practice, even though all studies for longitudinal count data with a fixed number of `failures' are conducted in a finite population, the NM model is typically preferred over the more accurate MIH model. Therefore, it is compelling to question how much information is lost when inferring based on the approximate NM model instead of the true MIH model. This loss is formally quantified by a measure known as deficiency. In this paper, asymptotic bounds for the deficiencies between the MIH and NM experiments are derived, as well as between the MIH and the corresponding multivariate normal experiments with identical mean-covariance structures, see Proposition~\ref{prop:3}. To achieve this goal, a local approximation for the log-ratio of the MIH and NM probability mass functions is proved first (Proposition~\ref{prop:1}) after which Hellinger distance bounds are established (Proposition~\ref{prop:2}). A similar study was done in \cite{MR4361955} for the multivariate hypergeometric distribution in relation with the multinomial distribution.

    \vspace{2mm}
    Here is an outline of the paper. In Section~\ref{sec:definitions}, the necessary definitions and notations are introduced. The results are stated in Section~\ref{sec:main} and proved in Section~\ref{sec:proofs}.

\section{Definitions}\label{sec:definitions}

    Define $\|\bb{x}\|_1 = \sum_{i=1}^d |x_i|$ to be the $\ell^1$ norm on $\R^d$. Following Equation~(2.7) of \citet{MR345316}, the probability mass function of the $\mathrm{MIH}\hspace{0.2mm}(N,n,\bb{p})$ distribution is supported on
    \begin{equation*}
        \mathbb{K}_d = \N_0^d \cap \textstyle\varprod_{i=1}^d [0, N p_i],
    \end{equation*}
    and defined as
    \begin{equation}\label{eq:inverse.hypergeometric.pmf}
        \begin{aligned}
            P_{N,n,\bb{p}}(\bb{k})
            &= \frac{(k_+ - 1)! (N - k_+)!}{N!} \frac{(N q)!}{(n - 1)! (N q - n)!} \prod_{i=1}^d \frac{(N p_i)!}{k_i! (N p_i - k_i)!} \\
            &= \frac{\Gamma(k_+)}{\Gamma(n) \prod_{i=1}^d k_i!} \frac{(N - k_+)!}{N!} \prod_{i=1}^{d+1} \frac{(N p_i)!}{(N p_i - k_i)!}, \quad \bb{k}\in \mathbb{K}_d,
        \end{aligned}
    \end{equation}
    where $k_{d+1} \equiv n \leq N q$ is fixed, $k_+ \equiv \sum_{i=1}^{d+1} k_i = n + \|\bb{k}\|_1$, and $\Gamma(x) = \int_0^{\infty} t^{x-1} e^{-t} \rd t, ~x\in (0,\infty)$, denotes Euler's gamma function. Following for example Chapter~36 of \citet{MR1429617}, the probability mass function of the $\mathrm{NM}\hspace{0.2mm}(n,\bb{p})$ distribution is defined by
    \begin{equation}\label{eq:negative.multinomial.pmf}
        Q_{n,\bb{p}}(\bb{k}) = \frac{\Gamma(k_+)}{\Gamma(n) \prod_{i=1}^d k_i!} \prod_{i=1}^{d+1} p_i^{k_i}, \quad \bb{k}\in \N_0^d,
    \end{equation}
    As mentioned above, this distribution represents the exact same as \eqref{eq:inverse.hypergeometric.pmf} above, except that the population from which the $n$ objects are drawn is infinite ($N = \infty$). For general information about these two distributions, refer to \cite{MR171341,MR345316,MR3559356,MR4358040,MR4448853}.

    \vspace{2mm}
    Throughout the paper, the notation $u = \OO(v)$ means that there exists a universal constant $C\in (0,\infty)$ such that $\limsup_{N\to \infty} |u / v| \leq C$. Whenever $C$ might depend on some parameter $\gamma$, say, one writes $u = \OO_{\gamma}(v)$.

\section{Results}\label{sec:main}

    The first main result is an asymptotic expansion for the log-ratio of the MIH probability mass function \eqref{eq:inverse.hypergeometric.pmf} over the corresponding negative multinomial probability mass function \eqref{eq:negative.multinomial.pmf}.

    \begin{proposition}[Local approximation]\label{prop:1}
        Assume that the integers $n,N\in \N$ satisfy $n\leq (1/2) \, N q$ and $\bb{p}\in N^{-1} \, \N^d$ holds, and pick any $\gamma\in (0,1)$.
        Then, uniformly for $\bb{k}\in \mathbb{K}_d$ such that
        \begin{equation*}
            \max_{i\in \{1,\dots,d\}} \frac{k_i}{p_i} \leq \gamma N,
        \end{equation*}
        one has, as $N\to \infty$,
        \begin{equation*}
            \begin{aligned}
                \log\left\{\frac{P_{N,n,\bb{p}}(\bb{k})}{Q_{n,\bb{p}}(\bb{k})}\right\}
                &= \frac{1}{N} \left[\left\{\frac{k_+^2}{2} - \frac{k_+}{2}\right\} - \sum_{i=1}^{d+1} \frac{1}{p_i} \left(\frac{k_i^2}{2} - \frac{k_i}{2}\right)\right] \\
                &~~+ \frac{1}{N^2} \left[\left\{\frac{k_+^3}{6} - \frac{k_+^2}{4} + \frac{k_+}{12}\right\} - \sum_{i=1}^{d+1} \frac{1}{p_i^2} \left(\frac{k_i^3}{6} - \frac{k_i^2}{4} + \frac{k_i}{12}\right)\right] \\
                &~~+ \OO_{\gamma}\left[\frac{1}{N^3} \left\{k_+^4 + \sum_{i=1}^{d+1} \left(\frac{k_i^4}{p_i^3} + \frac{k_i^2}{p_i^2}\right)\right\}\right].
            \end{aligned}
        \end{equation*}
    \end{proposition}

    The local approximation above together with the Hellinger distance bound in \cite{doi:10.1111/stan.12328} between jittered negative multinomials and the corresponding multivariate normals allow us to derive an upper bound on the Hellinger distance between the probability measure on $\R^d$ induced by a MIH random vector jittered by a uniform random vector on $(-1/2,1/2)^d$ and the probability measure on $\R^d$ induced by a multivariate normal random vector with the same mean and covariances as the negative multinomial distribution in \eqref{eq:negative.multinomial.pmf}.

    \begin{proposition}[Hellinger distance upper bounds]\label{prop:2}
        Assume that the integers $n,N\in \N$ satisfy $n\leq (1/2) \, N q$ and $\bb{p}\in N^{-1} \, \N^d$ holds.
        Let $\bb{K}\sim \mathrm{MIH}\hspace{0.2mm}(N,n,\bb{p})$, $\bb{L}\sim \mathrm{NM}\hspace{0.2mm}(n,\bb{p})$, and $\bb{U}, \bb{V}\sim \mathrm{Uniform}\hspace{0.2mm}(-1/2,1/2)^d$, where $\bb{K}$, $\bb{L}$, $\bb{U}$ and $\bb{V}$ are assumed to be jointly independent.
        Define $\bb{X} = \bb{K} + \bb{U}$ and $\bb{Y} = \bb{L} + \bb{V}$, and let $\widetilde{\PP}_{N,n,\bb{p}}$ and $\widetilde{\QQ}_{n,\bb{p}}$ be the laws of $\bb{X}$ and $\bb{Y}$, respectively.
        Also, let $\QQ_{n,\bb{p}}$ be the probability measure on $\R^d$ induced by the $\mathrm{Normal}_d(n \bb{p} / q, n \Sigma_{\bb{p}})$ distribution, where $\Sigma_{\bb{p}} = \mathrm{diag}(\bb{p} / q) - (\bb{p} / q) (\bb{p} / q)^{\top}$.
        Then, as $N\to \infty$, one has
        \begin{equation*}
            \mathcal{H}(\widetilde{\PP}_{N,n,\bb{p}},\widetilde{\QQ}_{n,\bb{p}})
            \leq \sqrt{200 \, d \exp\left\{- \frac{q \min(p_1,\ldots,p_d) N^2}{100 n}\right\} + \OO\left(\frac{d \, n^2}{N q^2}\right)},
        \end{equation*}
        and
        \begin{equation*}
            \mathcal{H}(\widetilde{\PP}_{N,n,\bb{p}},\QQ_{n,\bb{p}})
            \leq \mathcal{H}(\widetilde{\PP}_{N,n,\bb{p}},\widetilde{\QQ}_{n,\bb{p}}) + \OO\left\{\frac{d}{\sqrt{n \, q \, \min(p_1,\ldots,p_d)}}\right\},
        \end{equation*}
        where $\mathcal{H}(\cdot,\cdot)$ denotes the Hellinger distance.
    \end{proposition}

    \begin{remark}
        Considering the connection between the Hellinger distance and various other probability metrics as highlighted by \citet[p.~421]{doi:10.2307/1403865}, the above bounds remain applicable when substituting the Hellinger distance with any of these probability metrics: discrepancy metric, Kolmogorov metric, L\'evy metric, Prokhorov metric, total variation.
    \end{remark}

    Since the Markov kernel that jitters a random vector by a uniform on $(-1/2,1/2)^d$ is easily inverted (round off each component of the vector to the nearest integer), then one finds, as a consequence of the Hellinger distance bounds in Proposition~\ref{prop:2}, an upper bound on the maximum of the two deficiencies (namely $\delta(\mathscr{P}_b,\mathscr{Q}_b)$ and $\delta(\mathscr{Q}_b,\mathscr{P}_b)$) between MIH and multivariate normal experiments.

    \begin{proposition}[Deficiencies upper bounds]\label{prop:3}
        For any given $b\in (0,1)$, let
        \begin{equation*}
            \Theta_b = \left\{\bb{p}\in N^{-1} \, \N^d : \|\bb{p}\|_1 < 1 ~~ \text{and} ~~ \min(p_1,\ldots,p_d,q) \geq b\right\}.
        \end{equation*}
        Assume that the integers $n,N\in \N$ satisfy $n\leq (1/2) \, N b$.
        Define the experiments
        \begin{alignat*}{6}
            &\mathscr{P}_b
            &&= &&~\{\PP_{N,n,\bb{p}}\}_{\bb{p}\in \Theta_b}, \quad &&\PP_{N,n,\bb{p}} ~\text{is the measure induced by } \mathrm{MIH}\hspace{0.2mm}(N,n,\bb{p}), \\
            &\mathscr{Q}_b\hspace{-0.5mm}
            &&= &&~\{\QQ_{n,\bb{p}}\}_{\bb{p}\in \Theta_b}, \quad &&\QQ_{n,\bb{p}} ~\text{is the measure induced by } \mathrm{Normal}_d(n \bb{p} / q, n \Sigma_{\bb{p}}),
        \end{alignat*}
        where recall $\Sigma_{\bb{p}} = \mathrm{diag}(\bb{p} / q) - (\bb{p} / q) (\bb{p} / q)^{\top}$.
        Then, for $N\geq n^3 / d$, one has the following upper bounds on the deficiencies between the experiments $\mathscr{P}_b$ and $\mathscr{Q}_b$,
        \begin{equation*}
            \max\{\delta(\mathscr{P}_b,\mathscr{Q}_b),\delta(\mathscr{Q}_b,\mathscr{P}_b)\} \leq C_b \, \frac{d}{\sqrt{n}},
        \end{equation*}
        where $C_b$ is a positive real constant that depends only on $b$, and
        \begin{equation}\label{eq:def:deficiency.one.sided}
            \begin{aligned}
                \delta(\mathscr{P}_b,\mathscr{Q}_b)
                &= \inf_{T_1} \sup_{\bb{p}\in \Theta_b} \bigg\|\int_{\mathbb{K}_d} T_1(\bb{k}, \cdot \, ) \, \PP_{N,n,\bb{p}}(\rd \bb{k}) - \QQ_{n,\bb{p}}\bigg\|, \\
                \delta(\mathscr{Q}_b,\mathscr{P}_b)
                &= \inf_{T_2} \sup_{\bb{p}\in \Theta_b} \bigg\|\PP_{N,n,\bb{p}} - \int_{\R^d} T_2(\bb{z}, \cdot \, ) \, \QQ_{n,\bb{p}}(\rd \bb{z})\bigg\|.
            \end{aligned}
        \end{equation}
        In the above equation, the infima are taken over all Markov kernels $T_1 : \mathbb{K}_d \times \mathscr{B}(\R^d) \to [0,1]$ and $T_2 : \R^d \times \mathscr{B}(\mathbb{K}_d) \to [0,1]$.
    \end{proposition}

    Now, consider the following multivariate normal experiments with independent components
    \begin{alignat*}{6}
        &\overline{\mathscr{Q}}_b
        &&= &&~\{\overline{\QQ}_{n,\bb{p}}\}_{\bb{p}\in \Theta_b}, \quad &&\overline{\QQ}_{n,\bb{p}} ~\text{is the measure induced by } \mathrm{Normal}_d\{n \bb{p} / q, n \mathrm{diag}(\bb{p} / q)\}, \\
        &\mathscr{Q}_b^{\star}\hspace{-0.5mm}
        &&= &&~\{\QQ_{n,\bb{p}}^{\star}\}_{\bb{p}\in \Theta_b}, \quad &&\QQ_{n,\bb{p}}^{\star} ~\text{is the measure induced by } \mathrm{Normal}_d\{\hspace{-0.75mm}\sqrt{n \bb{p} / q}, \mathrm{diag}(\bb{1}/4)\},
    \end{alignat*}
    where $\bb{1} = (1,1,\dots,1)^{\top}$,
    then \citet[Section~7]{MR1922539} showed, using a variance stabilizing transformation, that
    \begin{equation}\label{eq:LeCam.distance.indep.normals}
        \max\{\delta(\mathscr{Q}_b,\overline{\mathscr{Q}}_b),\delta(\overline{\mathscr{Q}}_b,\mathscr{Q}_b)\} \leq C_b \, \sqrt{\frac{d}{n}} \quad \text{and} \quad \max\{\delta(\overline{\mathscr{Q}}_b,\mathscr{Q}_b^{\star}),\delta(\mathscr{Q}_b^{\star},\overline{\mathscr{Q}}_b)\} \leq C_b \, \frac{d}{\sqrt{n}},
    \end{equation}
    with proper adjustments to the definition of the deficiencies in \eqref{eq:def:deficiency.one.sided}.

    The corollary below presents deficiencies upper bounds between the MIH and the above multivariate normal experiments with independent components. The statement follows directly from Proposition~\ref{prop:3}, Equation~\eqref{eq:LeCam.distance.indep.normals} and the triangle inequality.

    \begin{corollary}\label{cor:4}
        With the same notation as in Proposition~\ref{prop:3}, one has, for $N\geq n^3 / d$,
        \begin{equation*}
            \max\{\delta(\mathscr{P}_b,\overline{\mathscr{Q}}_b),\delta(\overline{\mathscr{Q}}_b,\mathscr{P}_b)\} \leq C_b \, \frac{d}{\sqrt{n}} \quad \text{and} \quad \max\{\delta(\mathscr{P}_b,\mathscr{Q}_b^{\star}),\delta(\mathscr{Q}_b^{\star},\mathscr{P}_b)\} \leq C_b \, \frac{d}{\sqrt{n}},
        \end{equation*}
        where $C_b$ is a positive real constant that depends only on $b$.
    \end{corollary}

\section{Proofs}\label{sec:proofs}

    \begin{proof}[\bf Proof of Proposition~\ref{prop:1}]
        Throughout the proof, the parameter $n\in \N$ satisfies $n\leq \gamma N$ for some fixed $\gamma\in (0,1)$, and the asymptotic expressions are valid as $N\to \infty$.
        Also, assume that $\bb{p}\in N^{-1} \, \N^d$ and $\bb{k}\in \mathbb{K}_d$.
        Using Stirling's formula,
        \begin{equation*}
            \log m! = \frac{1}{2} \log(2\pi) + (m + \tfrac{1}{2}) \log m - m + \frac{1}{12 m} + \OO(m^{-3}), \quad m\to \infty,
        \end{equation*}
        see, e.g., \citet[p.257]{MR0167642}, and taking logarithms in \eqref{eq:inverse.hypergeometric.pmf} and \eqref{eq:negative.multinomial.pmf}, one obtains
        \begin{align*}
            \log\left\{\frac{P_{N,n,\bb{p}}(\bb{k})}{Q_{n,\bb{p}}(\bb{k})}\right\}
            &= \log (N - k_+)! - \log N! \\[-2mm]
            &\quad+ \sum_{i=1}^{d+1} \log (N p_i)! - \sum_{i=1}^{d+1} \log (N p_i - k_i)! - \sum_{i=1}^{d+1} k_i \log p_i \\[2mm]
            &= (N - k_+ + \tfrac{1}{2}) \log (N - k_+) - (N + \tfrac{1}{2}) \log N \\
            &\quad+ \sum_{i=1}^{d+1} (N p_i + \tfrac{1}{2}) \log (N p_i) - \sum_{i=1}^{d+1} (N p_i - k_i + \tfrac{1}{2}) \log (N p_i - k_i) - \sum_{i=1}^{d+1} k_i \log p_i \\
            &\quad+ \frac{1}{12 N} \left[\left(1 - \frac{k_+}{N}\right)^{-1} - 1 + \sum_{i=1}^{d+1} \frac{1}{p_i} \left\{1 - \left(1 - \frac{k_i}{N p_i}\right)^{-1}\right\}\right] \\
            &\quad+ \OO\left(\frac{1}{N^3} \left[\left(1 - \frac{k_+}{N}\right)^{-3} + 1 + \sum_{i=1}^{d+1} \frac{1}{p_i^3} \left\{1 + \left(1 - \frac{k_i}{N p_i}\right)^{-3}\right\}\right]\right).
        \end{align*}
        Simple algebraic manipulations yield
        \begin{align*}
            \log\left\{\frac{P_{N,n,\bb{p}}(\bb{k})}{Q_{n,\bb{p}}(\bb{k})}\right\}
            &= N \left(1 - \frac{k_+}{N}\right) \log \left(1 - \frac{k_+}{N}\right) + \frac{1}{2} \log \left(1 - \frac{k_+}{N}\right) \\
            &\quad- \sum_{i=1}^{d+1} N p_i \, \left(1 - \frac{k_i}{N p_i}\right) \log \left(1 - \frac{k_i}{N p_i}\right) - \frac{1}{2} \sum_{i=1}^{d+1} \log \left(1 - \frac{k_i}{N p_i}\right) \\
            &\quad+ \frac{1}{12 N} \left[\left(1 - \frac{k_+}{N}\right)^{-1} - 1 + \sum_{i=1}^{d+1} \frac{1}{p_i} \left\{1 - \left(1 - \frac{k_i}{N p_i}\right)^{-1}\right\}\right] \\
            &\quad+ \OO_{\gamma}\left\{\sum_{i=1}^{d+1} \frac{1}{(N p_i)^3}\right\}.
        \end{align*}
        By applying the following Taylor expansions, valid for $|x| \leq \eta < 1$,
        \begin{equation}\label{eq:prop:p.k.expansion.eq.Taylor}
            \begin{aligned}
                (1 - x) \log (1 - x) &= - x + \frac{x^2}{2} + \frac{x^3}{6} + \OO\left\{(1 - \eta)^{-3} |x|^4\right\}, \\[0mm]
                - \frac{1}{2} \log (1 - x) &= \frac{x}{2} + \frac{x^2}{4} + \OO\left\{(1 - \eta)^{-3} |x|^3\right\}, \\[2mm]
                (1 - x)^{-1} &= 1 + x + \OO\left\{(1 - \eta)^{-3} |x|^2\right\},
            \end{aligned}
        \end{equation}
        one has
        \begin{align*}
            \log\left\{\frac{P_{N,n,\bb{p}}(\bb{k})}{Q_{n,\bb{p}}(\bb{k})}\right\}
            &= N \left\{- \left(\frac{k_+}{N}\right) + \frac{1}{2} \left(\frac{k_+}{N}\right)^2 + \frac{1}{6} \left(\frac{k_+}{N}\right)^3 + \OO_{\gamma}\left(\bigg|\frac{k_+}{N}\bigg|^4\right)\right\} \\
            &\quad- \left\{\frac{1}{2}\left(\frac{k_+}{N}\right) + \frac{1}{4} \left(\frac{k_+}{N}\right)^2 + \OO_{\gamma}\left(\bigg|\frac{k_+}{N}\bigg|^3\right)\right\} \\
            &\quad+ \sum_{i=1}^{d+1} N p_i \, \left\{\left(\frac{k_i}{N p_i}\right) - \frac{1}{2} \left(\frac{k_i}{N p_i}\right)^2 - \frac{1}{6} \left(\frac{k_i}{N p_i}\right)^3 + \OO_{\gamma}\left(\bigg|\frac{k_i}{N p_i}\bigg|^4\right)\right\} \\
            &\quad+ \sum_{i=1}^{d+1} \left\{\frac{1}{2} \left(\frac{k_i}{N p_i}\right) + \frac{1}{4} \left(\frac{k_i}{N p_i}\right)^2 + \OO_{\gamma}\left(\bigg|\frac{k_i}{N p_i}\bigg|^3\right)\right\} \\
            &\quad+ \frac{1}{12 N} \left[\left\{\frac{k_+}{N} + \OO_{\gamma}\left(\bigg|\frac{k_+}{N}\bigg|^2\right)\right\} + \sum_{i=1}^{d+1} \frac{1}{p_i} \left\{- \frac{k_i}{N p_i} + \OO_{\gamma}\left(\bigg|\frac{k_i}{N p_i}\bigg|^2\right)\right\}\right] \notag \\
            &\quad+ \OO_{\gamma}\left\{\sum_{i=1}^{d+1} \frac{1}{(N p_i)^3}\right\}.
        \end{align*}
        After rearranging some terms, the conclusion follows.
    \end{proof}

    \begin{proof}[\bf Proof of Proposition~\ref{prop:2}]
        Define
        \begin{equation*}
            A_{N,n,\bb{p}}(\gamma) = \left\{\bb{k}\in \mathbb{K}_d : \max_{i\in \{1,\dots,d\}} \frac{k_i}{p_i} \leq \gamma N\right\} + \left(-\frac{1}{2}, \frac{1}{2}\right)^d, \quad \gamma\in (0,1).
        \end{equation*}
        By the comparison of the squared Hellinger distance with Kullback-Leibler divergence on page~726 of \citet{MR1922539}, one already knows that
        \begin{equation}\label{eq:first.bound.Hellinger.distance}
            \mathcal{H}^2(\widetilde{\PP}_{N,n,\bb{p}},\widetilde{\QQ}_{n,\bb{p}}\} \leq 2 \, \PP\left\{\bb{X}\in A_{N,n,\bb{p}}^c(1/2)\right\} + \EE\left[\log\left\{\frac{\rd \widetilde{\PP}_{N,n,\bb{p}}}{\rd \widetilde{\QQ}_{n,\bb{p}}}(\bb{X})\right\} \, \ind_{A_{N,n,\bb{p}}(1/2)}(\bb{X})\right].
        \end{equation}
        By applying a union bound together with a concentration bound for the (univariate) inverse hypergeometric distribution using the mean and variance derived on page 373 of in \citet{MR345316}, one gets, for $N$ large enough (recall that $n\leq (1/2) N q$),
        \begin{equation}\label{eq:concentration.bound}
            \begin{aligned}
                \PP\left\{\bb{X}\in A_{N,n,\bb{p}}^c(3/4)\right\}
                &\leq \sum_{i=1}^d \PP\left\{K_i > (3/4) N p_i\right\} \\
                &\leq \sum_{i=1}^d \PP\left(\frac{K_i - n p_i / q}{\sqrt{n (p_i / q) \{1 + (p_i / q)\}}} > \frac{(3/4) N p_i - n p_i / q}{\sqrt{n (p_i / q) \{1 + (p_i / q)\}}}\right) \\
                &\leq \sum_{i=1}^d \PP\left(\frac{K_i - n p_i / q}{\sqrt{n (p_i / q) \{1 + (p_i / q)\}}} > \frac{(1/4) N p_i}{\sqrt{n (p_i / q)}}\right) \\
                &\leq 100 \, d \exp\left\{- \frac{q \min(p_1,\ldots,p_d) N^2}{100 n}\right\}.
            \end{aligned}
        \end{equation}
        To estimate the expectation in \eqref{eq:first.bound.Hellinger.distance}, note that if $P_{N,n,\bb{p}}(\bb{x})$ and $Q_{n,\bb{p}}(\bb{x})$ denote the density functions associated with $\widetilde{\PP}_{N,n,\bb{p}}$ and $\widetilde{\QQ}_{n,\bb{p}}$ (i.e., $P_{N,n,\bb{p}}(\bb{x})$ is equal to $P_{N,n,\bb{p}}(\bb{k})$ whenever $\bb{k}\in \mathbb{K}_d$ is closest to $\bb{x}$, and analogously for $Q_{n,\bb{p}}(\bb{x})$), then, for $N$ large enough, one has
        \begin{equation*}
            \begin{aligned}
                \left|\EE\left[\log\left\{\frac{\rd \widetilde{\PP}_{N,n,\bb{p}}}{\rd \widetilde{\QQ}_{n,\bb{p}}}(\bb{X})\right\} \, \ind_{A_{N,n,\bb{p}}(1/2)}(\bb{X})\right]\right|
                &= \left|\EE\left[\log\left(\frac{P_{N,n,\bb{p}}(\bb{X})}{Q_{n,\bb{p}}(\bb{X})}\right) \, \ind_{A_{N,n,\bb{p}}(1/2)}(\bb{X})\right]\right| \\
                &\leq \EE\left[\bigg|\log\left(\frac{P_{N,n,\bb{p}}(\bb{K})}{Q_{n,\bb{p}}(\bb{K})}\right)\bigg| \, \ind_{A_{N,n,\bb{p}}(3/4)}(\bb{X})\right].
            \end{aligned}
        \end{equation*}
        By Proposition~\ref{prop:1} with $\gamma = 3/4$, one finds
        \begin{equation*}
            \begin{aligned}
                \EE\left[\log\left\{\frac{\rd \widetilde{\PP}_{N,n,\bb{p}}}{\rd \widetilde{\QQ}_{n,\bb{p}}}(\bb{X})\right\} \, \ind_{\{\bb{X}\in A_{N,n,\bb{p}}(3/4)\}}\right]
                &= \OO\left\{\frac{n^2 + \sum_{i=1}^d \EE(K_i^2)}{N} + \sum_{i=1}^{d+1} \frac{\EE(K_i^2)}{N p_i}\right\} \\
                &= \OO\left(\frac{n^2 + \sum_{i=1}^d n^2 p_i q^{-2}}{N} + \sum_{i=1}^{d+1} \frac{n^2 p_i q^{-2}}{N p_i}\right) \\
                &= \OO\left(\frac{d \, n^2}{N q^2}\right).
            \end{aligned}
        \end{equation*}
        Together with the concentration bound in \eqref{eq:concentration.bound}, one deduces from \eqref{eq:first.bound.Hellinger.distance} that
        \begin{equation*}
            \mathcal{H}^2(\widetilde{\PP}_{N,n,\bb{p}},\widetilde{\QQ}_{n,\bb{p}}) \leq 200 \, d \exp\left\{- \frac{q \min(p_1,\ldots,p_d) N^2}{100 n}\right\} + \OO\left(\frac{d \, n^2}{N q^2}\right).
        \end{equation*}
        which proves the first claim of the proposition.
        Furthermore, by Proposition~2 in \cite{doi:10.1111/stan.12328}, it is already established that
        \vspace{-1mm}
        \begin{equation*}
            \mathcal{H}(\widetilde{\QQ}_{n,\bb{p}},\QQ_{n,\bb{p}})
            = \OO\left\{\frac{d}{\sqrt{n \, q \, \min(p_1,\ldots,p_d)}}\right\},
        \end{equation*}
        which proves the second claim of the proposition.
        This ends the proof.
    \end{proof}

    \begin{proof}[\bf Proof of Proposition~\ref{prop:3}]
        By Proposition~\ref{prop:2} with the assumption $N\geq n^3 / d$, one obtains the desired bound on $\delta(\mathscr{P}_b,\mathscr{Q}_b)$ by choosing the Markov kernel $T_1^{\star}$ that adds the uniform jittering $\bb{U}\sim \mathrm{Uniform}\hspace{0.2mm}(-1/2,1/2)^d$ to $\bb{K}\sim \mathrm{MIH}\hspace{0.2mm}(N,n,\bb{p})$, namely
        \begin{equation*}
            \begin{aligned}
                T_1^{\star}(\bb{k},B) = \int_{(-\frac{1}{2},\frac{1}{2})^d} \ind_{B}(\bb{k} + \bb{u}) \rd \bb{u}, \quad \bb{k}\in \mathbb{K}_d, ~B\in \mathscr{B}(\R^d).
            \end{aligned}
        \end{equation*}
        To obtain the bound on $\delta(\mathscr{Q}_b,\mathscr{P}_b)$, it suffices to consider a Markov kernel $T_2^{\star}$ that inverts the effect of $T_1^{\star}$, i.e., rounding off every components of $\bb{Z}\sim \mathrm{Normal}_d(n \bb{p} / q, n \Sigma_{\bb{p}})$ to the nearest integer.
        Then, as explained by \citet[Section~5]{MR1922539}, one gets
        \begin{equation*}
            \begin{aligned}
                \delta(\mathscr{Q}_b,\mathscr{P}_b)
                &\leq \bigg\|\PP_{N,n,\bb{p}} - \int_{\R^d} T_2^{\star}(\bb{z}, \cdot \, ) \, \QQ_{n,\bb{p}}(\rd \bb{z})\bigg\| \\[0.5mm]
                &= \bigg\|\int_{\R^d} T_2^{\star}(\bb{z}, \cdot \, ) \int_{\mathbb{K}_d} T_1^{\star}(\bb{k}, \rd \bb{z}) \, \PP_{N,n,\bb{p}}(\rd \bb{k}) - \int_{\R^d} T_2^{\star}(\bb{z}, \cdot \, ) \, \QQ_{n,\bb{p}}(\rd \bb{z})\bigg\| \\
                &\leq \bigg\|\int_{\mathbb{K}_d} T_1^{\star}(\bb{k}, \cdot \, ) \, \PP_{N,n,\bb{p}}(\rd \bb{k}) - \QQ_{n,\bb{p}}\bigg\|,
            \end{aligned}
        \end{equation*}
        and one finds the same bound as on $\delta(\mathscr{P}_b,\mathscr{Q}_b)$ by Proposition~\ref{prop:2}.
    \end{proof}

\section*{Competing interests}

The author declares no competing interests.

\section*{Acknowledgments}

We express our gratitude to the referee for carefully reviewing our manuscript and providing valuable suggestions that have enhanced the presentation of this paper.

\section*{Funding}

F.\ Ouimet is supported by a CRM-Simons postdoctoral fellowship from the Centre de recherches math\'ematiques (Montr\'eal, Canada) and the Simons Foundation.

%
%

\phantomsection
\addcontentsline{toc}{chapter}{References}

\bibliographystyle{authordate1}
\bibliography{Ouimet_2023_LLT_inverse_hypergeometric_bib}

\end{document}